%% file: root.tex
\documentclass{ifacconf}

\makeatletter
\let\old@ssect\@ssect % Store how ifacconf defines \@ssect
\makeatother
\usepackage{natbib}   
\usepackage[hidelinks]{hyperref}
\usepackage{orcidlink}
\makeatletter
\def\@ssect#1#2#3#4#5#6{%
  \NR@gettitle{#6}% Insert key \nameref title grab
  \old@ssect{#1}{#2}{#3}{#4}{#5}{#6}% Restore ifacconf's \@ssect
}
\makeatother
\usepackage{graphicx}      % include this line if your document contains figures
\usepackage{lineno}
\usepackage{amsmath,amssymb,amsfonts} %for equations
\usepackage{siunitx}
\usepackage{booktabs}
\usepackage{xcolor}

\usepackage{fancyhdr}
\usepackage[export]{adjustbox}
\newcommand{\V}[1]{\boldsymbol{#1}}%use this command for vectors
\newcommand{\M}[1]{\mathbf{#1}}%use this command for matrices
\let\oldalign\align
\let\oldendalign\endalign

\renewenvironment{align}
  {\linenomathNonumbers\oldalign}
  {\oldendalign\endlinenomath}
\begin{document}

\begin{frontmatter}

\title{Optimal Control for Wind Turbine Wake Mixing on Floating Platforms} 

\thanks[footnoteinfo]{This work is part of the research programme ``Robust closed-loop wake steering for large densely spaced wind farms'' with project number 17512, which is (partly) financed by the Dutch Research Council (NWO).
This project is part of the Floatech project. The research presented in this paper has received funding from the European Union's Horizon 2020 research and innovation programme under grant agreement No. 101007142.
}

\author[First]{Maarten J. van den Broek\,\orcidlink{0000-0001-7396-4001}}
\author[First]{Daniel van den Berg\,\orcidlink{0000-0002-1623-1482}}
% \author{\dots} 
\author[Second]{Benjamin Sanderse\,\orcidlink{0000-0001-9483-1988}} 
\author[First]{Jan-Willem van Wingerden\,\orcidlink{0000-0003-3061-7442}}

\address[First]{
Delft Centre for Systems and Control, TU Delft, Delft, The~Netherlands (e-mail: m.j.vandenbroek@tudelft.nl)
}
\address[Second]{Scientific Computing, CWI, Amsterdam, The~Netherlands
}

\input{sections/abstract}


\end{frontmatter}
\thispagestyle{fancy}
%===============================================================================
% \linenumbers
\input{sections/introduction}
% \thispagestyle{plain}
\input{sections/methods}

\input{sections/results}

\input{sections/conclusion}

\input{sections/notes}

% \begin{ack}
% Place acknowledgments here.
% \end{ack}
% \bibliography{ifacconf}             % bib file to produce the bibliography
\bibliographystyle{ifacconf}
\bibliography{./wind-farm-control.bib,./references.bib}
                                                     % with bibtex (preferred)

\appendix
% \section{A summary of Latin grammar}    % Each appendix must have a short title.
% \section{Some Latin vocabulary}              % Sections and subsections are supported  
                                                                         % in the appendices.
\end{document}

%% file: sections/abstract.tex
%!TEX root=../main.tex

\begin{abstract}                % Abstract of not more than 250 words.
% Reduce aerodynamic wake interaction for increased wind farm power production.
%
Dynamic induction control is a wind farm flow control strategy that utilises wind turbine thrust variations to accelerate breakdown of the aerodynamic wake and improve downstream turbine performance.
%
%Floating wind turbines introduce additional dynamics and challenges compared to conventional fixed-bottom wind turbines.
However, when floating wind turbines are considered, additional dynamics and challenges appear that make optimal control difficult.
%
%This work builds on previous results using adjoint optimisation with a free-vortex wake model for economic model-predictive control for dynamic wind farm flow control solutions. %to optimise dynamic induction control for fixed-bottom turbines to improve mean power production.
%
%In this work, the aerodynamic model is coupled with linear floating platform hydrodynamics with pitch and surge degrees of freedom.
%In this work, we propose an adjoint optimisation framework for non-linear economic model-predictive control by coupling the aerodynamic wake model to floating platform hydrodynamics.
In this work, we propose an adjoint optimisation framework for non-linear economic model-predictive control, which utilises a novel coupling of an existing aerodynamic wake model to floating platform hydrodynamics.
% word? implements? utilises?
%
% Implementation with automatic differentiation allows increased flexibility in model development and enables gradient-based optimisation.
 % for application in an economic model-predictive control setting to maximise mean power production in a two-turbine case.
%
Analysis of the frequency response for the coupled model shows that it is possible to achieve wind turbine thrust variations without inducing large motion of the rotor.
Using economic model-predictive control, we find dynamic induction results that lead to an improvement of \textcolor{black}{\SI{7}{\%}} %\SI{4.8}{\%}
over static induction control, where
the dynamic controller stimulates wake breakdown with only small variations in rotor displacement.
%
%This novel model formulation couples wind turbine wake aerodynamics and floating platform hydrodynamics to explore dynamic wind farm flow control strategies on floating wind turbines.
%The dynamic control results demonstrate the possibility for induction control on floating turbines and provide a starting point for studies with higher fidelity models.
%
This novel model formulation provides a starting point for the adaptation of dynamic wind farm flow control strategies for floating wind turbines.
\end{abstract}

%% file: sections/introduction.tex
%!TEX root=../root.tex

\section{Introduction}\label{sec:introduction}

Offshore parcels suitable for fixed-bottom wind turbines are limited, leading to wind turbines being built in large, densely-spaced wind farms~\citep{VanWingerden2020a}.
Floating wind turbines provide an opportunity to extend the suitable space for wind farm construction beyond shallow waters (i.e. deeper than 50 metres).
%However, floating platforms introduce new challenges in terms of wind turbine and farm-scale control.
% To enable growth of offshore wind, transition to floating wind turbines.

% Why do we do dynamic induction control?
% What is dynamic induction control?
Wind farm flow control aims to reduce the negative effects of aerodynamic interaction between wind turbines in a farm with strategies such as wake redirection by yaw misalignment, induction control, and wake mixing strategies~\citep{Meyers2022}.
These control strategies have mostly been developed on fixed-bottom turbines.
% , such as yaw-roll lock \textcolor{blue}{[cite]} and platform pitch instabilities \textcolor{blue}{[cite]}.
% These unique dynamics 
% require specific attention in 
The implementation of wind farm flow control strategies for floating wind turbines requires specific attention as floating platforms introduce additional dynamics, instabilities, and resonant modes.
% Mostly developed on fixed bottom turbines.
% Unique dynamics for floating platform wind turbines deserve specific attention
% yaw roll-lock, pitch instability with regular wind turbine controller.
% Why consider floating turbines?
% What are floating turbines?
% How do floating turbines behave?

% Once category of wind farm control strategies is dynamic induction control, where Control strategies for wind farm control can be roughly divided in three categories: w
The current work focuses on dynamic induction control - the application of wind turbine thrust variations, and consequent induction variations, to stimulate breakdown of the aerodynamic wake behind the turbine.
% The origin of these dynamic induction control strategies lies with optimal control studies with adjoint optimisation of large-eddy simulations~\citep{Goit2015}.
% In a step towards practical application, the control signals found through optimisation have been reduced to sinusoidal thrust signals applied using collective pitch control. 
Optimal control studies with adjoint optimisation of large-eddy simulations~\citep{Goit2015} provided the basis for dynamic induction control signals, which, in a step towards practical application, have been reduced to sinusoidal thrust signals applied using collective pitch control. 
These have been found to improve wake recovery in a study with large-eddy simulations~\citep{Munters2018a} and in wind tunnel experiments~\citep{Frederik2020a}.  

% The patterns found through optimal control studies with adjoint large-eddy simulations provided the basis for these dynamic induction control methods, although these simulations are too computationally expensive for real-time control applications~\citep{Munters2018,Munters2018a}.

Recently, adjoint optimisation with a free-vortex wake representation of a wind turbine wake using a two-dimensional actuator-disc model has been used to find dynamic induction control signals in a non-linear economic model-predictive control setting \citep{Broek2022arxiv}.
This study provided a cumbersome manual derivation of the discrete adjoint to enable gradient-based optimisation. 
It is computationally inexpensive compared to optimisation using large-eddy simulations.
% Dynamic induction control signals to stimulate wake breakdown in a two-turbine case study.

% Model-based control optimsiation and adaptation to floating platforms?
The modelling of wake aerodynamics using free-vortex methods is suitable for adaptation to floating wind turbine platforms because the simulation uses Lagrangian particles, therefore allowing rotor motions without being limited by grid resolution and localised refinements.
The wakes of floating wind turbines have previously been modelled using a free-vortex ring method~\citep{Dong2019} and free-vortex wake models have been used to study the wake dynamics~\citep{Lee2019}.
However, these studies focused on the performance of a single turbine instead of considering aerodynamic wake interaction.
% Within the field of wind energy, free-vortex wake models have been used to study floating wind turbines and wake dynamics~\citep{Lee2019}.
% To study dynamic wake control methods and analyse wake stability~\citep{Brown2021}.

To explore the possibilities for dynamic induction control on floating wind turbines, we extend the free-vortex wake model by \citet{Broek2022arxiv} with a floating platform model and utilise automatic differentiation to replace manual derivation of the discrete adjoint.
The coupled aerodynamic and hydrodynamic model can then be used to explore optimal control of dynamic induction with gradient-based optimisation in an economic model-predictive control framework.

Contributions of this paper are:
(i) extension of the free-vortex wake model with floating platform dynamics,
(ii) implementation of the optimal control algorithm with automatic differentiation to avoid manual derivation of the adjoint, and
(iii) exploration of wake mixing using optimal control for wind turbines on floating platforms.
%in a comparison between fixed-bottom and floating platforms.

The remainder of this paper is structured as follows.
Section~\ref{sec:methods} describes the wake model for the floating wind turbine and methods for control optimisation.
Section~\ref{sec:results} discusses the results on wake mixing for floating platforms.
Finally, conclusions are presented in Section~\ref{sec:conclusion}.

%% file: sections/methods.tex
%!TEX root=../root.tex

% \section{Methods}\label{sec:methods}
\section{Modelling and Control}\label{sec:methods}
This section introduces the novel coupling of the aerodynamic wake model and hydrodynamic platform model for dynamic control optimisation.
Section~\ref{subsec:wake_model} describes, in brief, the use of the free-vortex wake to construct a 2D actuator-disc model for the representation of the wind turbine wake. %, adapted from \citet{Broek2022arxiv}.
The linear hydrodynamic model of a floating platform is presented in Section~\ref{subsec:tower_model} and the connection to the aerodynamic model is described in Section~\ref{subsec:modelcoupling}.
An objective function for mean power maximisation and receding horizon controller is defined in Section~\ref{subsec:control}, such that optimal control signals may be explored in this novel framework.

\subsection{Wind Turbine Wake Model}\label{subsec:wake_model}
The free-vortex method is used to to model the aerodynamic wake using  a two-dimensional (2D) actuator-disc representation of a wind turbine. The basis of this model is presented in \citet{Broek2022arxiv}.
The vorticity formulation requires the assumption of inviscid and incompressible flow, although diffusion may be approximated.
The actuator disc is assumed to be uniformly loaded so it only releases vorticity along its edge~\citep{Katz2001}.

\begin{figure}[!b]
	\centering
	\includegraphics[width=\linewidth]{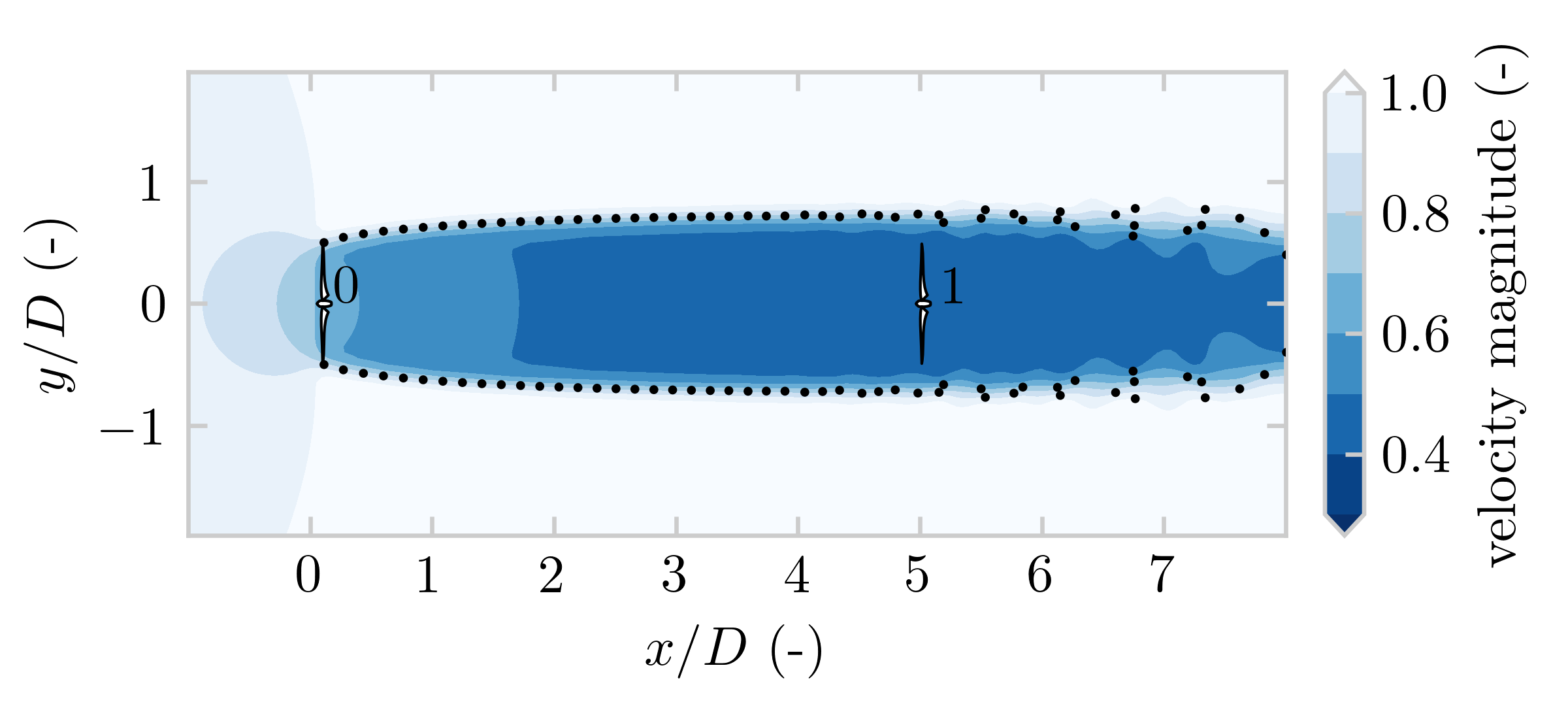}
	\caption{
	Illustration of a flow field at hub height with the wind turbine wake for turbine~0, where turbine~1 at $x/D=5$ downstream is modelled using the effective flow velocity.
	Units have been non-dimensionalised by rotor diameter and inflow velocity magnitude.
	The vortex points used for the velocity calculation are marked with black dots.
	% Mirroring of the vortex points and related circulation is used to model the effect of the ground.
	}\label{fig:wake_model}
\end{figure}

% The 2D wake model in this work represents a vertical slice ($xz$-plane), whereas \textcolor{blue}{[preprint]} considered a horizontal slice ($xy$-plane) at hub height.
% This new orientation allows inclusion of platform pitch in the wake aerodynamics.
% The ground effect -- no flow may pass through the ground or water surface -- is included by mirroring the system of vortex points \textcolor{blue}{[references]}.
% Figure~\ref{fig:wake_model} shows an example flow field and the vortex points used to calculate the wake effects of the wind turbine.
The 2D wake model in this work represents a horizontal slice ($xy$-plane) of the flow field at hub height.
At every time step, pairs of vortex points are released from the edge of the rotor that travel downstream as Lagrangian markers.
These vortex points are used to calculate the velocity deficit in the wake as illustrated in Figure~\ref{fig:wake_model}.
% \begin{itemize}
% 	\item 2D Actuator-disc model using the free-vortex Wake
% 	\item mirroring for ground effect
% 	\item set of vortex points at every time-step
% \end{itemize}

% Relevant assumptions:
% \begin{itemize}
% 	\item inviscid/potential flow
% 	\item uniformly loaded actuator disc{}
% 	\item lack of turbulence (diffusion may be approximated)
% \end{itemize}

The wake dynamics are modelled as a non-linear state-space system in discrete time for the state update at a time-step $k$,
\begin{align}
	\V{q}_{k+1} &= f(\V{q}_k, \V{m}_k)\,, \label{eq:state_update}\\
		\V{y}_{k} &= g(\V{q}_k, \V{m}_k)\,, \label{eq:state_output}
\end{align}
with the state vector $\V{q}\in\mathbb{R}^{n_\mathrm{s}}$, the control vector $\V{m}\in\mathbb{R}^{n_\mathrm{c}}$, and the output vector $\V{y}\in\mathbb{R}^{n_\mathrm{t}}$.
The total number of states to describe the wake is $n_\mathrm{s}$, the total number of controls is $n_\mathrm{c}$ and the number of turbines
% both with a wake model and virtual, 
for which power output is estimated
is $n_\mathrm{t}$.

The velocity $\V{u}_\mathrm{i}\in\mathbb{R}^{2}$ induced at point $\V{x}_0\in\mathbb{R}^{2}$ by a single vortex element located at $\V{x}_1\in\mathbb{R}^{2}$,  with vortex strength $\Gamma$, is calculated with the Biot-Savart law as 
\begin{align}
	\V{u}_\mathrm{i}(\V{x}_0, \V{x}_1) = 
\begin{bmatrix}-r_y\\r_x\end{bmatrix}
\left(
\frac{\Gamma}{2\pi}
\frac{1}{||\V{r}||^2}
\right)
\left(
1 - \exp\left(-\frac{||\V{r}||^2}{\sigma^2}\right)
\right)\,,
\label{eq:ui_2d}
\end{align}
where the relative position $\V{r}$ is
\begin{align}
	\V{r} = \V{x}_1 - \V{x}_0\,.
\end{align}
A Gaussian core with core size $\sigma$ is included to regularise singular behaviour of the induced velocity close to the vortex element.
The total velocity induced by all vortex elements is calculated by summation of all individual contributions.
At every time step, a new set of vortex points is initialised at the edge of the rotor and all other points are propagated downstream.
% For a more complete definition of the aerodynamic model, the reader is referred to \citet{Broek2022arxiv}.

The thrust $T$ is calculated according to the mean velocity at the rotor plane $\V{u}_\mathrm{r}$,
\begin{align}
	T = c_\mathrm{t}'(a)\cdot \frac{1}{2}\rho A_\mathrm{r}(\V{n}\cdot\V{u}_\mathrm{r})^2\,,
\end{align}
with local thrust coefficient $c_\mathrm{t}'$, air density $\rho$, rotor swept area $A_\mathrm{r}$ and the rotor normal vector $\V{n}$.
% include \V{n}(\phi)??
Similarly, power $P$ is calculated with the local power coefficient $c_\mathrm{p}'$ as
\begin{align}
	P = c_\mathrm{p}'(a)\cdot \frac{1}{2}\rho A_\mathrm{r}(\V{n}\cdot\V{u}_\mathrm{r})^3\,.
\end{align}
The thrust and power coefficient are a function of the induction factor $a$ based on momentum theory,
\begin{align}
	c_\mathrm{t}'(a) &= \left\{\begin{array}{ll}
	\frac{4a(1-a)}{(1-a)^2} = \frac{4a}{1-a} 	& \text{if } a\leq a_\mathrm{t}\,,\\
	\frac{c_\mathrm{t1}-4(\sqrt{c_\mathrm{t1}}-1)(1-a)}{(1-a)^2} 
	& \text{if } a>a_\mathrm{t}\,,
	\end{array}
	\right.\\
	c_\mathrm{p}'(a) &= \frac{4a}{1-a}\,,
\end{align}
with the parameter $c_\mathrm{t1}=2.3$.
The thrust coefficient calculation is based on momentum theory with a transition to a linear approximation for high induction values which is an empirical correction based on the Glauert correction~\citep{Burton2001}.

For the performance estimate of virtual turbines, based on the velocity estimate without including their impact on the flow, the rotor velocity is lowered using the induction factor,
\begin{align}
	\V{u}_\mathrm{r}^*=(1-a)\V{u}_\mathrm{r} \,.
\end{align}

% \textcolor{blue}{Important frequency at $St=0.2-0.3$, given wind speed and rotor diameter, corresponds to SI{0.01}{Hz}.
% Extend frequency discussion and notes on choice of time step}

\begin{table}[!b]
	\centering
	\caption{Parameters for wind turbine, floating platform, and numerical configuration.}\label{tab:physical_param}
	\begin{tabular}{lcr}
	\toprule
		Rotor diameter 	& $D$ & \SI{178.3}{m}\\
		Nacelle height 	& $h$ & \SI{119.}{m}\\
		% Rotor tilt 			& $\theta$ & \SI{5.}{\degree}\\
		Total mass 			& $m_0$ & \SI{1.1e6}{kg}\\
		Mass moment of inertia & $I_{yy,0}$ & \SI{3.9e10}{kg.m^2}\\
		Added mass 			& $\Delta m$ & \SI{2.8e7}{kg}\\
		Added inertia 			& $\Delta I_{yy}$ & \SI{1.1e10}{kg.m^2}\\
		\midrule
		Pitch stiffness &$k_{\phi}$& \SI{6.2e9}{N.m.rad^{-1}} \\
		Pitch damping &$c_{\phi}$& \SI{7.3e8}{N.m.s.rad^{-1}}\\
		Surge stiffness &$k_{x}$& \SI{8.3e4}{N.m^{-1}} \\
		Surge damping &$c_x$& \SI{1.7e5}{N.s.m^{-1}}\\
		\midrule
% 	\end{tabular}
% \end{table}
% % wind turbine mass 1 106 954 kg
% % platform total mass 28 828 000 kg
% % platform
% 	% surge 0.005 Hz, heave 0.06 Hz. pitch 0.04 Hz.  tt fore-aft 0.4 Hz
% % linear damping 
% % surge 1.7 e5 Ns/m
% % heave 1.35 e6 Ns/m
% % pitch 7.3 e8 Nms 
% % yaw 1.09 e8 Nms
% % added mass wrt SWL
% % surge 2.75 e7 kg
% % pitch 1.085 e10 kg m2
% % surge stiffness -> mooring lines?
% % stiffness about SWL initial position / @ x = 19.57 m
% % pitch - 55
% 	% -2.0002e8 / -1.675e8 Nm/rad
% % surge - 11
%  % 8.3283e4 / 7.6888e4 N/m ?
%  \begin{table}[!hbt]
% \centering
% \caption{Parameters for numerical simulation}\label{tab:sim_param}
% \begin{tabular}{lcr}
% \toprule
Time step - floater & $\Delta t_\mathrm{f}$ & \SI{0.9}{s}\\
Time step - wake & $\Delta t_\mathrm{w}$ & \SI{3.6}{s}\\
Number of rings & $n_\mathrm{r}$ & 60\\
Vortex core size & $\sigma$ & \SI{17.8}{m}\\
Air density & $\rho$ &\SI{1.225}{kg\cdot m^3} \\
Inflow velocity & $\V{u}_{\infty}$ & \SI{10.0}{m.s^{-1}}\\
\bottomrule
\end{tabular}
\end{table}

\subsection{Floating Platform Model}\label{subsec:tower_model}
\begin{figure}[!b]
	\centering
	\includegraphics[width=\linewidth]{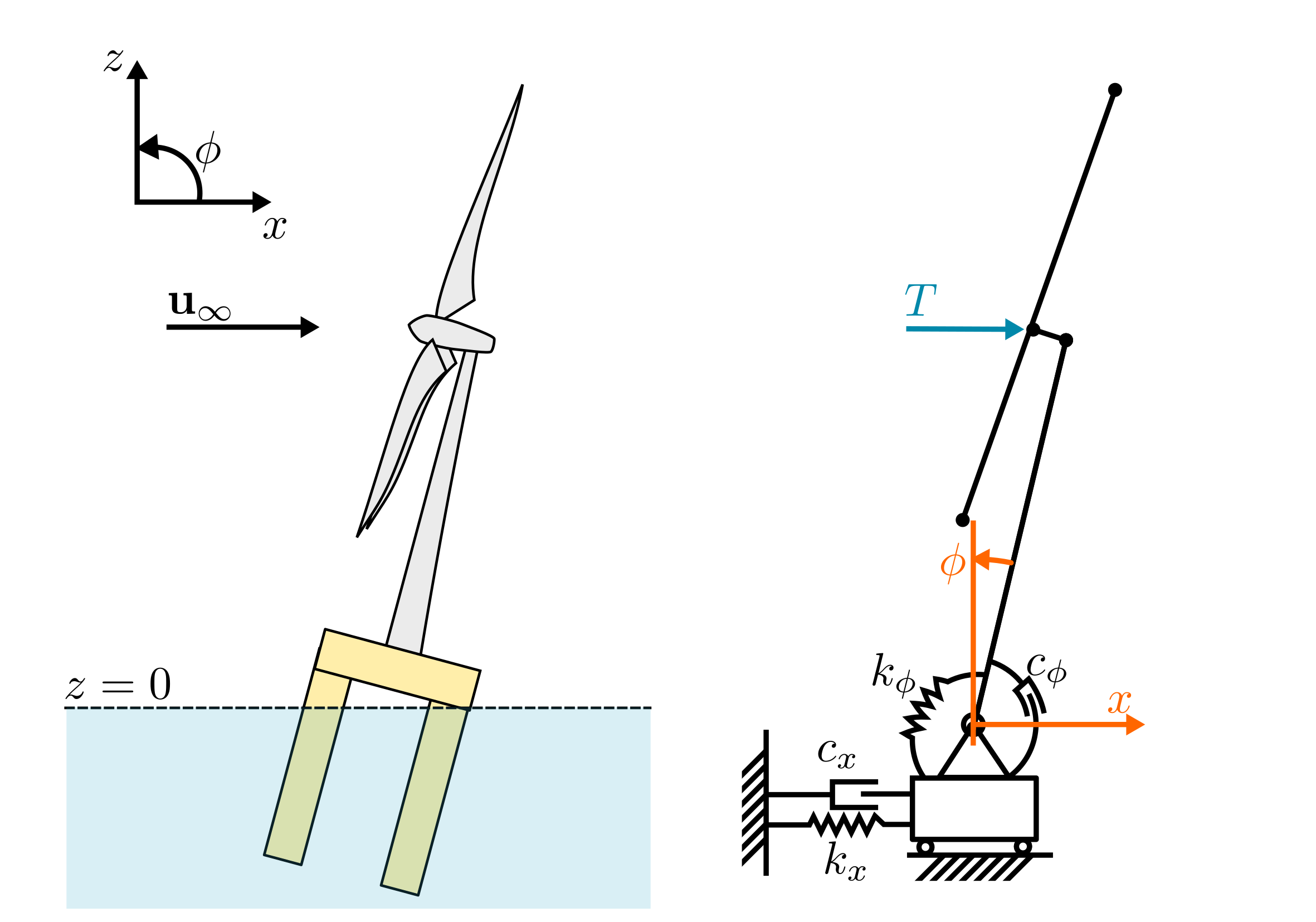}
	\caption{
	Diagram of the simplified model for the dynamics of the floating wind turbine platform.
	Platform surge and pitch dynamics are modelled using, respectively, a translational and rotational mass-spring-damper system.
	}\label{fig:FBD}
\end{figure}

The displacement of the rotor due to the dynamics of a floating base is calculated from platform motions in pitch (rotation along $y$) and surge (translation in $x$) degrees of freedom.
These are the main motions influencing wind turbine thrust and power \citep{Lee2019}.
The heave motion (translation in $z$) of the platform is neglected.
The hydrodynamics of the platform motions are modelled using two decoupled mass-spring-damper systems, translational and rotational for surge and pitch, respectively.{}
Figure~\ref{fig:FBD} illustrates the model for the floating platform dynamics.
The parameters used in this study are based on a triple-spar platform with a \SI{10}{MW} reference turbine \citep{Lemmer2018} and are listed in Table~\ref{tab:physical_param}.

The assumption of linear platform hydrodynamics may be strong, but by neglecting non-linear hydrodynamics, we can put emphasis on the complexity of wake aerodynamics.
%non-linear hydrodynamics are outside of the scope of this study, which emphasises the complexities of wake aerodynamics.
Additionally, the tower and blades are assumed to be rigid.
Decoupling of pitch and surge relies on assuming that translation happens at the centre of rotation. 
% \textcolor{blue}{[valid?]}
The tilt and pitch angles are assumed small, and therefore the thrust is assumed aligned with the $x$-axis.

% Assumptions:
% \begin{itemize}
% 	\item Linear second-order dynamics (hydrodynamics are non-linear, but necessary for simplicity)
% 	\item Neglect heave motion
% 	\item Decoupled response between pitch and surge DOF
% 	\item thrust orthogonal to rotor / orthogonal to tower??
% 	\item rigid tower
% 	\item small angles?
% \end{itemize}

The equation of motion for platform pitch follows from a balance of moments around the centre of rotation,
\begin{align}
	I_\mathrm{yy}\ddot{\phi} &= -Th-k_{\phi}\phi -c_{\phi}\dot{\phi}\,,
\end{align}
with the mass moment of inertia $I_{yy}$, pitch angle~$\phi$, thrust~$T$, tower height~$h$, rotational hydrodynamic stiffness and damping $k_{\phi}$ and $c_{\phi}$, respectively.
The mass moment of inertia is the sum of physical inertia $I_{yy,0}$ and a hydrodynamic added inertia term $\Delta I_{yy}$.

The equation of motion for surge follows from the balance of forces in $x$-direction,
\begin{align}
	m\ddot{x} &= T - k_{x} x - c_x \dot{x}\,,
\end{align}
with total mass $m$, displacement of the centre of rotation $x$, and translational hydrodynamic stiffness and damping $k_x$ and $c_x$, respectively.
The total mass is the sum of physical mass $m_{0}$ and a hydrodynamic added mass term $\Delta{m}$.

Both equations of motion are put into a linear state-space form $\dot{\V{x}}=\M{A}\V{x}+\M{B}\V{u}$ as
\begin{align}
		\begin{bmatrix}
		\dot{\phi} \\ \ddot{\phi}
	\end{bmatrix} &= 
	\begin{bmatrix}
	0	& 1\\
	-k_{\phi}/I_{yy}	&  -c_{\phi}/I_{yy}\\
	\end{bmatrix} \begin{bmatrix}
		\phi \\ \dot{\phi}
	\end{bmatrix}+\begin{bmatrix}
		0 \\ -h/I_{yy}\\ 
	\end{bmatrix}T\,,\\
			\begin{bmatrix}
		\dot{x} \\ \ddot{x}
	\end{bmatrix} &= 
	\begin{bmatrix}
	0	& 1\\
	-k_{x}/m	&  -c_{x}/m\\
	\end{bmatrix} \begin{bmatrix}
		x \\ \dot{x}
	\end{bmatrix}+\begin{bmatrix}
		0 \\ 1/m\\ 
	\end{bmatrix}T\,.
\end{align}

These are converted to discrete-time state-space form using a zero-order hold with a time-step $\Delta t_\mathrm{f}$ for application in the numerical model framework according to
% $\M{A}_\mathrm{d}=\exp(\M{A}\Delta t)$ and $\M{B}_{\mathrm{d}}=\M{A}^{-1}(\M{A}_\mathrm{d}-I)\M{B}$,
\begin{align}
  \M{A}_\mathrm{d}&=\exp(\M{A}\Delta t_\mathrm{f}) \,,\\
  \M{B}_{\mathrm{d}}&=\M{A}^{-1}(\M{A}_\mathrm{d}-\M{I})\M{B}\,,\\
	\V{x}_{k+1} &= \M{A}_\mathrm{d}\V{x}_{k} + \M{B}_\mathrm{d}T\,. \label{eq:platform_state}
\end{align}

\subsection{Model Coupling}\label{subsec:modelcoupling}
There exists a two-way coupling between the two models that is evaluated at every time step in the wake simulation.
The thrust force is calculated from the difference between the effective wind speed and the motion of the rotor as calculated by the floating platform model. 
The position of the rotor from the floating platform model also specifies the position where the vortex points are released in the wake model.

The sampling time needs to be chosen such that dynamics at the natural frequency of the system are adequately captured.
The resonant frequencies for the rigid floating platform are $\omega_{\phi}=\textcolor{black}{\SI{0.056}{Hz}}$ in the pitch and $\omega_x=\textcolor{black}{\SI{0.0085}{Hz}}$ in the surge degree of freedom.
The sampling frequency of the wake model is insufficient to properly capture the floating dynamics for platform pitching motion.
This is resolved by running the floater dynamics at a higher frequency, i.e. at a smaller time step $\Delta t_\mathrm{w}=4\Delta t_\mathrm{f}$.
The wake model is assumed to be constant during these smaller time steps for updating the thrust input to the floating platform dynamics.

Running the wake model at a higher sampling frequency is undesired because that would require equally more vortex points and more steps in the optimisation horizon.
That quickly makes optimisation slow and cumbersome.

% \textcolor{blue}{Note on difference with sampling time for $St=0.2-0.3$ from dynamic induction control.}

% High sampling frequencies and associated small time steps are problematic for optimisation over long horizons.
% The memory cost and compilation time for the AD gradient calculations make optimisation prohibitive -- same goes for experimentation on a laptop.
% From the frequency sweep, one can conclude that the pitch frequency dynamics are likely not to be dominant for power performance estimates.
% A coarse time resolution that captures the main aerodynamic properties may then be a good choice for control optimisation.

% For the given platform, natural frequencies are approximately: pitch at \SI{0.056}{Hz} and surge at \SI{0.0085}{Hz}.
% To represent pitch resonance, sampling frequency is minimally \SI{0.1}{Hz}, or $f=1.8$ assuming $D=\SI{180}{m}$ and $u=\SI{10}{m.s^{-1}}$.
% The associated time-step is $\SI{0.55}>0.20$, so pitch response should be captured in sampling time as used in [preprint].
% At nine steps per period, the signal is a bit rough.

% We may adjust for different dynamical time-scales by subsampling the simulations.
% Tower dynamics are cheap to compute, so can easily be run at twice the sampling frequency of the aerodynamic model.
% Assume the wake model is constant inbetween time-steps for updating thrust and power as inputs to tower dynamics.

\subsection{Control Optimisation}\label{subsec:control}
This work follows a non-linear economic model-predictive control (EMPC) approach for maximisation of mean power production as introduced in \citet{Broek2022arxiv} and briefly summarised here.
Receding horizon control is implemented without terminal constraints by implementing a prediction horizon that is sufficiently long to remove the impact of finite horizon effects on the implemented control solution \citep{Grune2013}.

The non-linear EMPC formulation requires optimisation of an objective function at every time step.
We construct a scalar objective function $J:\mathbb{R}^{n_\mathrm{s}n_\mathrm{c}}\rightarrow\mathbb{R}$, that combines the output power for both turbines in $\V{y}_k$ and the change in control signal $\V{\Delta m}_k = \V{m}_k-\V{m}_{k-1}$ as
\begin{align}
	 %J(\V{q}_k,\V{m}_k) = \M{Q}\V{y}_k(\V{q}_k,\V{m}_k) + \V{\Delta m}_k^{\operatorname{T}}\M{R}\V{\Delta m}_k \,,
	 J(\V{q}_k,\V{m}_k) = \M{Q}\V{y}_k + \V{\Delta m}_k^{\operatorname{T}}\M{R}\V{\Delta m}_k \,,
\end{align}
with output weight $\M{Q}\in\mathbb{R}^{1\times n_\mathrm{t}}$ and input weight $\M{R}\in\mathbb{R}^{n_\mathrm{c}\times n_\mathrm{c}}$.
The output weight is chosen element-wise negative ($\M{Q}<0$) such that minimisation of the objective maximises mean power production and the input weight is chosen positive ($\M{R}>0$) to penalise actuation cost and smoothen the optimisation landscape.
We then construct the unconstrained optimisation problem, subject to wake and platform dynamics, over a horizon of $N_\mathrm{h}$ steps, starting from step $k_0$,
\begin{align}
\label{eq:objective}
  \min_{\V{m}_k} \sum_{k=k_0}^{k_0+N_\mathrm{h}} J(\V{q}_k,\V{m}_k) \,,
\end{align}
% subject to the wake dynamics \eqref{eq:state_update}, \eqref{eq:state_output} and platform dynamics \eqref{eq:platform_state}, 
to find the optimal controls $\V{m}_k$ for $k\in[k_0;k_0+N_\mathrm{h}]$.

	% \V{q}_{k+1} &= f(\V{q}_k, \V{m}_k)\,, \label{eq:state_update}\\
	% 	\V{y}_{k} &= g(\V{q}_k, \V{m}_k)\,, \label{eq:state_output}

% The full definition of the finite-horizon optimal control problem is omitted here for reasons of space, but can be reconstructed from the optimisation problem presented in \citet{Broek2022arxiv}.
The optimisation problem is solved by using a gradient-based approach with the Adam optimiser using the default parameters \citep{Kingma2015}.
The current work replaces the manual derivation of the adjoint equations by utilising reverse mode automatic differentiation with Zygote~\citep{Zygote} to construct the gradient of the objective function.
% The gradient of the objective function is evaluated with reverse mode automatic differentiation using Zygote \textcolor{blue}{[ref]}.
% This is in contrast to the previous work in \textcolor{blue}{[preprint]} that utilised a manual derivation of the adjoint equations.

%% file: sections/results.tex
%!TEX root=../root.tex
\section{Results and Discussion}\label{sec:results}
This section presents the optimal control signals and characterisation of system dynamics for the novel coupled model framework and receding horizon control strategy described in the previous section.
A characterisation of the frequency response of the floating platform and wake aerodynamics is provided in Section~\ref{subsec:freqresp}.
This is followed by a discussion of the results from the receding horizon optimisation of controls in Section~\ref{subsec:empc_results}.

\subsection{Frequency Response}\label{subsec:freqresp}
The frequency domain characteristics of the coupled wake and floating platform model are estimated by actuating the turbine induction with a chirp signal from \textcolor{black}{\SI{0.0025}{Hz} to \SI{0.1}{Hz}} over \SI{30000}{s}, with a mean value of $\bar{a}=0.28$ and amplitude of $a=0.05$.
The FFT of the signals is used to estimate the frequency response functions.
The estimation is performed for a floating platform with both pitch and surge degrees of freedom and for %pitch and surge individually.
a single degree-of-freedom platform with either pitch or surge motion.

\begin{figure}[!b]
\centering
\includegraphics[width=\linewidth]{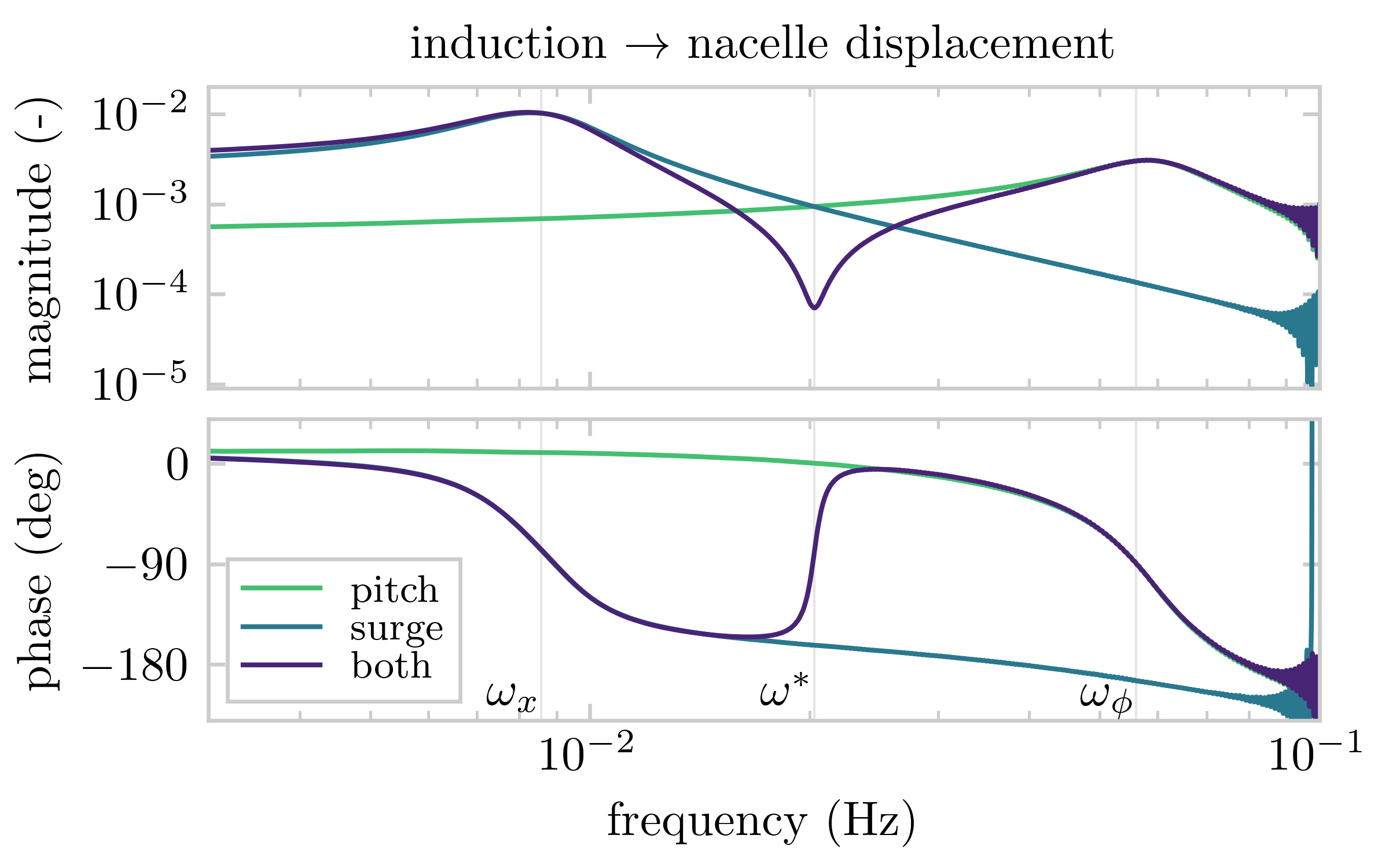}
\caption{
   Frequency response from induction excitation to tower top motion for different degrees of freedom of the floating platform. %, estimated using a chirp signal from \textcolor{black}{\SI{0.001}{Hz} to \SI{0.1}{Hz}}.
   Resonance in the platform motion can be excited for both surge and pitch degrees of freedom resulting in large displacements at the rotor.
   However, the anti-resonant peak shows the two motions can combine to lead to a relatively stationary rotor while maintaining thrust variations.
}\label{fig:towertop_tf}
\end{figure}

The response of wind turbine motion in response to dynamic induction variations is characterised by the response in nacelle displacement in Figure~\ref{fig:towertop_tf}.
The resonance in surge and pitch modes of the platform motion appears as peaks in the magnitude of the tower top response.
Actuation at these frequencies results in large motions of the rotor.
There is a cross-over around \textcolor{black}{$\omega^*=\SI{0.02}{Hz}$} where the combination of counter-phase pitch and surge motion results in a strongly reduced transmission to tower top displacement.
This indicates the possibility for dynamic induction actuation on floating wind turbines with minimal movement of the rotor.

\begin{figure}[!b]
\centering
\includegraphics[width=\linewidth]{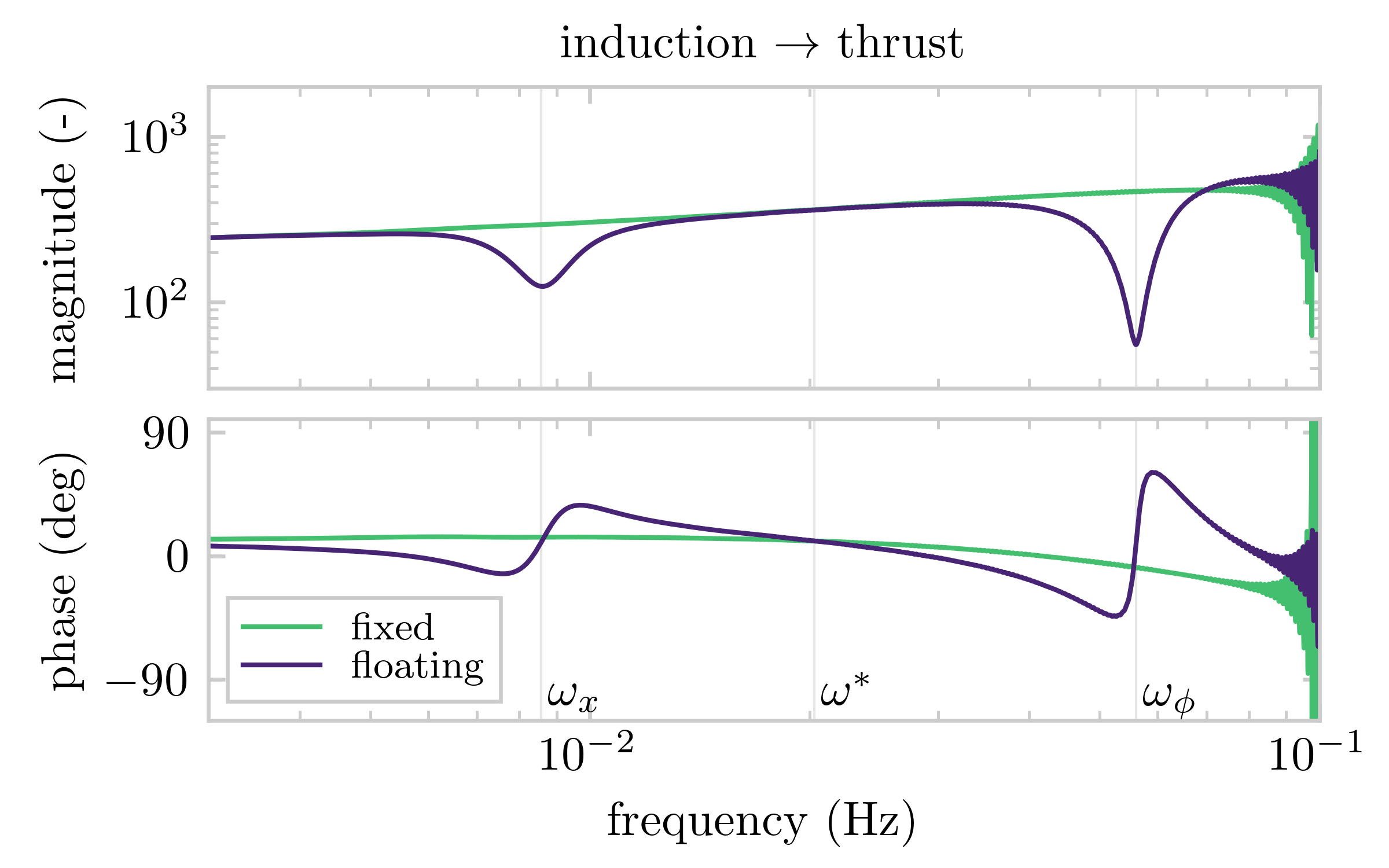}
\caption{
   Frequency response from induction excitation to turbine thrust for a bottom-fixed and floating wind turbine. %, estimated using a chirp signal from \textcolor{black}{\SI{0.001}{Hz} to \SI{0.1}{Hz}}.
   Resonance in the platform motion for surge and pitch reduces the magnitude of thrust variations in response to induction input.
}\label{fig:thrust_tf}
\end{figure}

The effects of platform resonance also appear in the frequency response of turbine thrust to induction excitation as illustrated in Figure~\ref{fig:thrust_tf}, when compared to a bottom-fixed turbine.
The motion of the rotor on the floating platform reduces the effective thrust variation that can be achieved as resonant platform modes are excited.
This reduction in thrust variation in response to the induction input is an important characteristic for finding optimal dynamic induction control signals on floating wind turbines.
Reduced thrust variations consequently reduce the ability to stimulate breakdown of the wake.

\begin{figure}[!b]
\centering
\includegraphics[width=\linewidth]{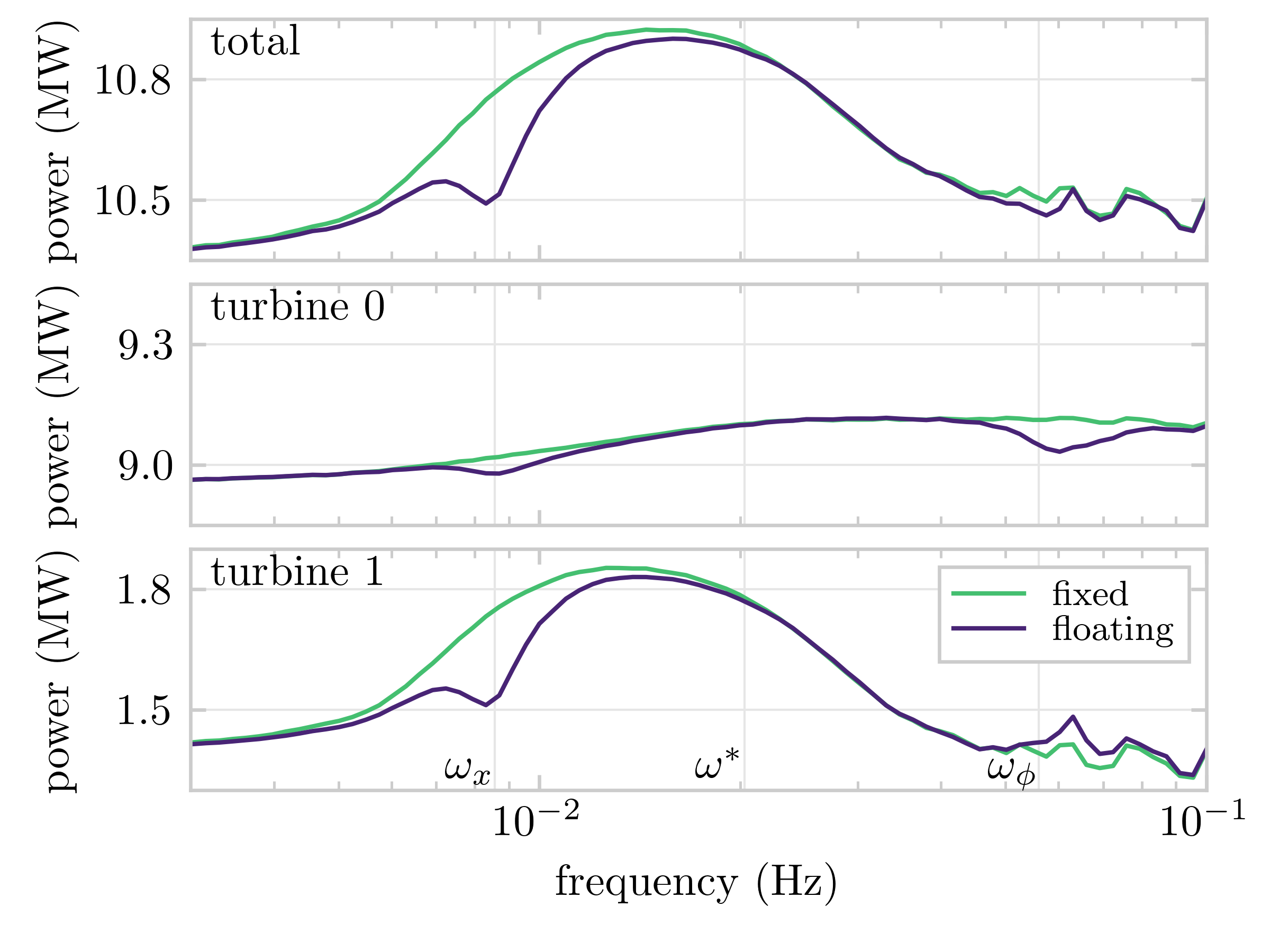}
\caption{
Mean power production in quasi-steady state for dynamic induction excitation over frequency,
comparing a simulation with a bottom-fixed turbine and a floating platform with pitch and surge degrees of freedom.
Maximum power is produced at an actuation frequency of \textcolor{black}{\SI{0.014}{Hz}, which equals a Strouhal number of $\mathrm{St}=0.26$}.
}\label{fig:maximalpowerfreq}  
\end{figure}

We evaluate the effect of induction variations 
on turbine 0, with turbine 1 at its individual optimum,
by taking a sinusoidal control signal that is slightly under-inductive,
\begin{align}
   a(t) = 0.28 + 0.05\sin (2\pi f t)\,,
\end{align}
and varying the frequency.
The mean power is calculated once the simulation has converged to a quasi-steady state for every frequency.
This stepped-frequency signal differs from the chirp signal in which this quasi-steady state for individual frequencies is not achieved.
The resulting power characteristic is illustrated in Figure~\ref{fig:maximalpowerfreq} for both fixed bottom turbines and for turbines on floating platforms.

For the illustrated sinusoidal sweep, the maximum mean power achieved is \SI{10.92}{MW} for the fixed turbines and \SI{10.90}{MW} for the turbines on floating platforms, at a frequency of \SI{0.014}{Hz} and \SI{0.016}{Hz}, respectively.
These frequencies correspond to a dimensionless frequency of $\mathrm{St}=0.26$ and $\mathrm{St}=0.28$, which is close to the frequency found in previous results \citep{Munters2018a,Broek2022arxiv}.
Dynamic induction actuation with a sinusoidal signal exceeds the maximum power production of \SI{10.44}{MW} that could otherwise be achieved in this model with a static induction signal for both fixed and floating turbines. 

The additional dynamics introduced by the floating platform appear slightly detrimental to the power that can be gained through thrust variations.
The excitation of platform natural frequencies reduces the thrust variations that can be achieved and consequently the ability to excite the aerodynamic breakdown of the wake on floating platforms.
Thrust variations can be achieved for frequencies in between the platform modes, especially those where the rotor is mostly stationary, i.e. the nacelle displacement due to induction variation is minimal.

\subsection{Economic Model-Predictive Control}\label{subsec:empc_results}
The EMPC case study is run over 300 time steps, or \SI{1080}{s}.
As in \citet{Broek2022arxiv}, the optimisation problem for control of turbine 0 is solved over a prediction horizon of \textcolor{black}{$N_\mathrm{h}=100$} steps, equivalent to \SI{360}{s}, and \textcolor{black}{50} optimiser iterations per step in the receding horizon.
% with a prediction horizon of \textcolor{black}{$N_\mathrm{h}=100$} steps, equivalent to \SI{360}{s}, and \textcolor{black}{50} optimiser iterations per step in the receding horizon as in~\citet{Broek2022arxiv}.
The objective weights 
$\M{Q}=\begin{bmatrix}
   -1 & -1 \\
\end{bmatrix}\,\si{MW^{-1}}$ and $\M{R}=\begin{bmatrix}
   \num{4.7e-2}
\end{bmatrix}$  are scaled to the physical problem. % for optimisation of the induction signal on turbine 0, assuming turbine 1 at $a=0.33$.
% The objective weights are set to
% $\M{Q}=\begin{bmatrix}
%    -1 & -1 \\
% \end{bmatrix}\,\si{MW^{-1}}$ and $\M{R}=\begin{bmatrix}
%    \num{4.7e-2}
% \end{bmatrix}$.
% As in \citet{Broek2022arxiv}, the optimisation problem is solved over a prediction horizon of \textcolor{black}{$N_\mathrm{h}=100$} steps, equivalent to \SI{360}{s}, and \textcolor{black}{50} optimiser iterations per step in the receding horizon.
% The control signal for turbine 0 is optimised with turbine 1 at its greedy optimum with $a=0.33$.
The first sample of the control solution is implemented as the simulation is advanced one time step, after which the control signal is shifted in time and re-optimised.
% Turbine 0 is controlled and turbine 1 operates at its greedy optimum.
% \textcolor{black}{The input weight $\M{R}=30$ is chosen higher than in the original work}

\begin{figure}[!b]
\centering
\includegraphics[width=\linewidth]{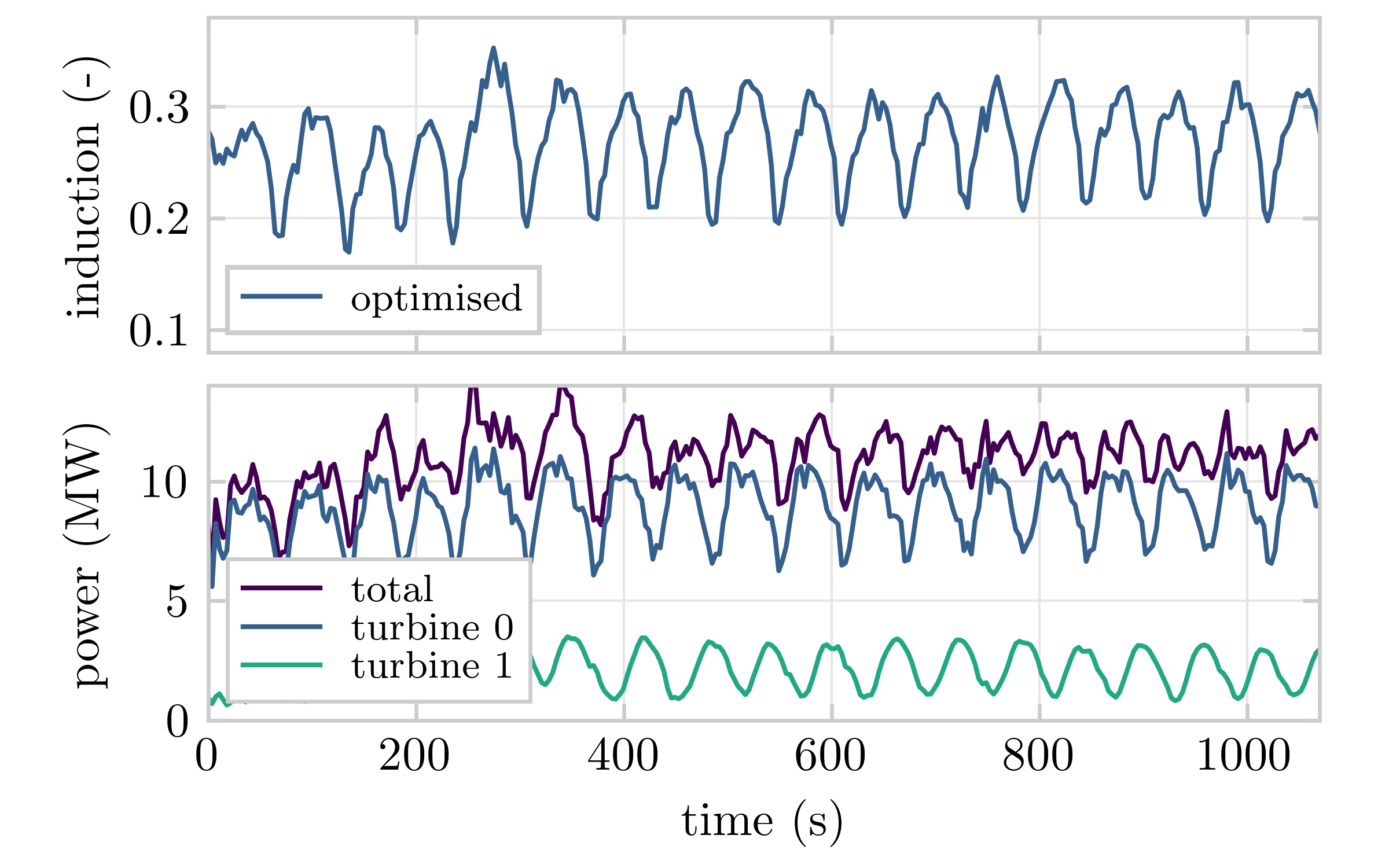}
\caption{
   The roughly periodic induction signal for turbine 0 is found through non-linear EMPC for maximisation of mean power production.
   The mean total power of \textcolor{black}{\SI{11.2}{MW}} towards the end of the simulation  is an improvement of \textcolor{black}{\SI{2.4}{\%}} over the maximum with sinusoidal actuation shown in Figure~\ref{fig:maximalpowerfreq}.
}\label{fig:empc_induction}
\end{figure}

\begin{figure}[!b]
\centering
\includegraphics[width=\linewidth]{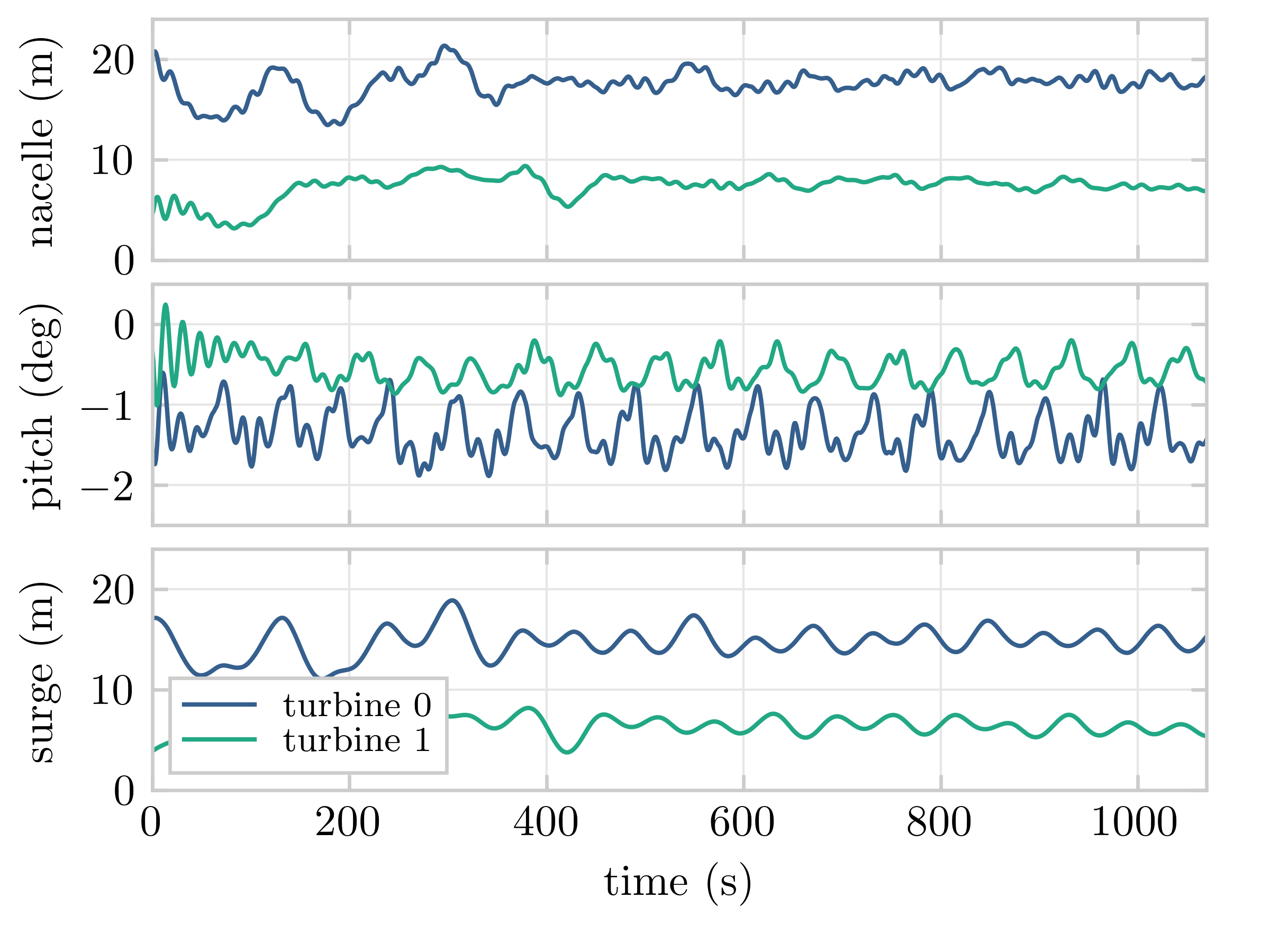}
\caption{The platform response associated with the control signal in Figure~\ref{fig:empc_induction}.
Dynamic nacelle displacement is minimal and mostly driven by platform surge motion. 
}\label{fig:empc_platform}
\end{figure}

The dynamic induction control signal and power estimate found through the receding horizon optimisation of controls are illustrated in Figure~\ref{fig:empc_induction}, with the associated platform motions and nacelle displacement in Figure~\ref{fig:empc_platform}.
The control signal becomes roughly periodic after the transients have passed 
and the power estimate converges to a mean total power of \textcolor{black}{\SI{11.2}{MW}} towards the end of the simulation.
This is an improvement of \textcolor{black}{\SI{6.9}{\%}} over what can be achieved with static induction control and \textcolor{black}{\SI{2.4}{\%}} over the maximum with sinusoidal actuation shown in Figure~\ref{fig:maximalpowerfreq}. 
The associated nacelle displacement is relatively small, although the dynamic variation of induction, related to thrust, does trigger a response in the platform motions. %pitch and surge motion.

The frequency content associated with the induction signal is presented in Figure~\ref{fig:empc_induction}.
% and compared with the baseline signal from \textcolor{black}{[preprint?]}.
The dominant frequency is \textcolor{black}{\SI{0.017}{Hz} ($\mathrm{St}=0.30$)}, which is close to the anti-resonance in the frequency response of the nacelle displacement (Figure~\ref{fig:towertop_tf}).
% or other baseline? 
\begin{figure}[!b]
\centering
\includegraphics[width=\linewidth]{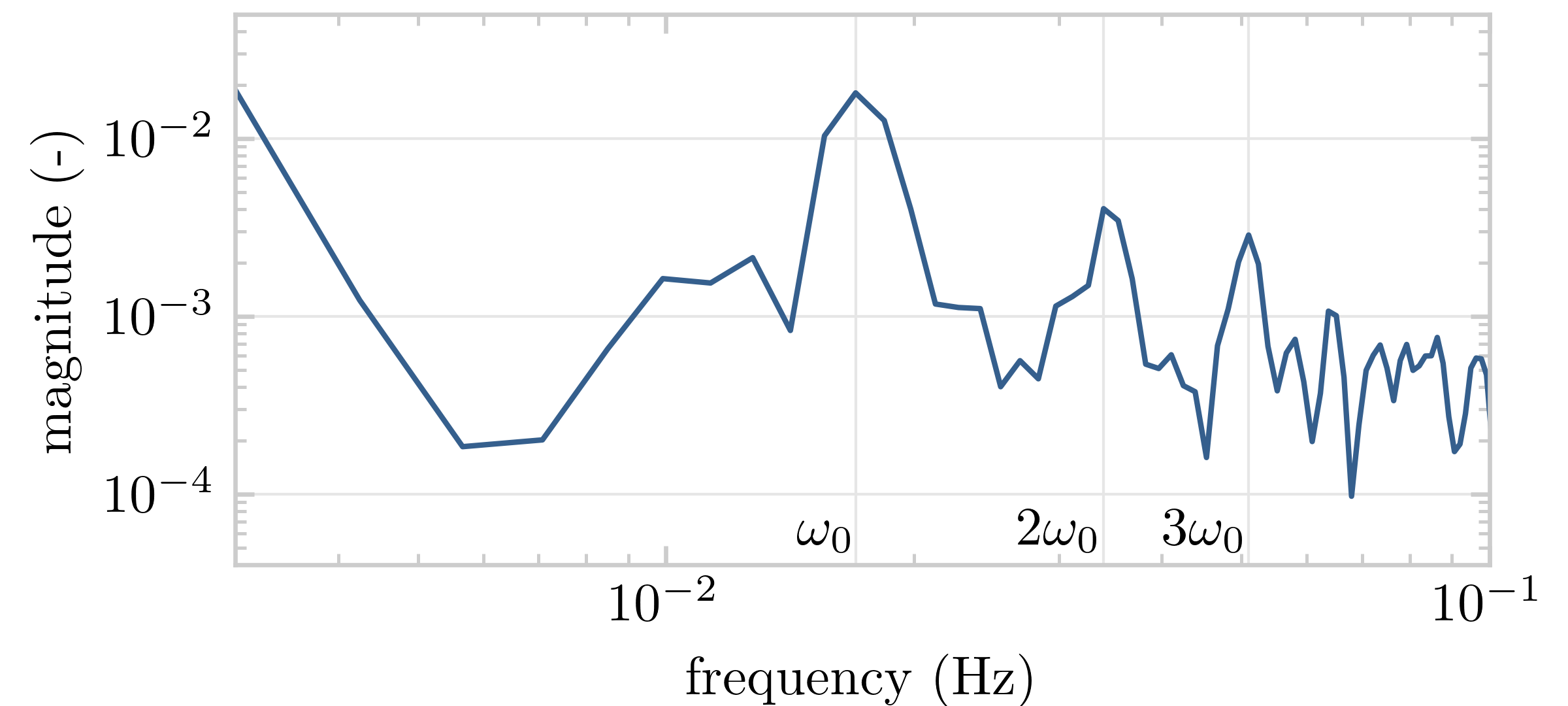}
\caption{The frequency content for the induction signal from Figure~\ref{fig:empc_induction} has a peak at $\omega_0=\SI{0.017}{Hz}$ and several harmonics.
The dominant peak is close to the frequency of \SI{0.02}{Hz} where transmission to tower top motion is minimal.
% \textcolor{black}{This signal is similar to that found for a fixed bottom turbine in the same configuration.}
}\label{fig:freq_content}
\end{figure}
Dynamic induction excitation at this anti-resonant frequency allows large thrust variations to promote the aerodynamic breakdown of the wake while minimising rotor movement.
The low sensitivity of upstream power to the precise frequency of thrust variation allows the optimisation space to find the optimal excitation frequency.
% Upstream power appears to have low sensitivity to the precise frequency of thrust variation. 
% This allows the aerodynamic optimum to appear in the frequency sweep and the optimisation.

%This excitation happens at a significantly higher frequency than the $St=0.25$ found by \citet{Munters2018a} for fixed-bottom turbines and which is commonly used in literature.
This excitation happens at significantly higher frequency than the $\mathrm{St}=0.20$ found by \citet{Broek2022arxiv} with the same model formulation, but for fixed-bottom turbines.
This indicates the significant impact of the floating platform on the wind farm flow control strategy.
It is important to note that the optimal frequency for aerodynamic excitation of wake breakdown should not be too close to the resonant frequencies for platform motion. 
The large displacement of the rotor would then lead to large power losses at the upstream turbine and limit the possibilities for stimulation of wake mixing due to loss of thrust variation.

%% file: sections/conclusion.tex
\section{Conclusion}\label{sec:conclusion}
% \textcolor{gray}{
% A conclusion section is not required. Although a conclusion may review
% the main points of the paper, do not replicate the abstract as the
% conclusion. A conclusion might elaborate on the importance of the work
% or suggest applications and extensions.
% }
% Contributions of this paper are
% (i) extension of the free-vortex wake model with floating platform dynamics,
% (ii) implementation with automatic differentiation to avoid manual derivation of the adjoint , and
% (iii) exploration of wake mixing using optimal control in a comparison between fixed-bottom and floating platforms.
The work presented in this paper provides a novel hybrid model-based approach for the optimisation of dynamic induction control signals for wind turbines on floating platforms to stimulate wake breakdown and improve mean power production under waked conditions.

A 2D free-vortex wake representation of the wind turbine wake using an actuator-disc model is coupled with a floating platform model with pitch and surge degrees of freedom.
% The implementation with automatic differentiation avoids manual adjoint derivations and enables flexibility in model development.
%
The exploration of the frequency response of the coupled wake and floating platform model shows that platform motions may lead to reduced thrust variations for dynamic induction control.
On the other hand, platform motions may also combine to allow large thrust variations with minimal rotor movement.
The control signal found in the non-linear EMPC demonstrates dynamic induction control for wind turbines on a floating platform, effectively stimulating aerodynamic wake breakdown without inducing large displacements of the rotor.

% \begin{itemize}
%    \item Adaptation of free-vortex wake model to $xz$-plane and inclusion of floating platform dynamics.
%    \item Implementation with automatic differentiation to avoid manual adjoint derivations and allow flexibility in model development.
%    \item Showcase example of dynamic induction control on floating platform. Results show that large thrust variations may stimulate wake breakdown without inducing large displacements of the rotor.
% \end{itemize}

Future work may explore validation with higher fidelity models or a co-design approach to simultaneously adapt platform parameters and control signals to improve dynamic control performance.

%% file: sections/notes.tex
%!TEX root=../main.tex

% \section*{Notes}